\documentclass[11pt]{amsart}
\usepackage{a4wide}
\usepackage{amsmath}
\usepackage[utf8]{inputenc}
\usepackage{amssymb}
\usepackage{amsopn}
\usepackage{epsfig}
\usepackage{amsfonts}
\usepackage{latexsym}
\usepackage{graphicx}
\usepackage{enumerate}
\usepackage[dvipsnames]{xcolor}
\usepackage{mathrsfs }
\usepackage{tikz}
\usetikzlibrary{arrows,automata}

\newcommand{\set}[1]{\left\{ #1 \right\}}

\newcommand{\R}{\mathbb{R}}

\newcommand{\N}{\mathbb{N}}

\newcommand{\f}{\infty}

\newcommand{\wt}[1]{\widetilde{#1}}
\newcommand{\ol}[1]{\overline{#1}}
\newcommand{\ul}[1]{\underline{#1}}
\newcommand{\ep}{\varepsilon}

\newcommand{\diam}{\mathrm{diam}\;}

\usepackage{amsthm}
\newtheorem{theorem}{Theorem}[section]
\newtheorem{proposition}[theorem]{Proposition}
\newtheorem{lemma}[theorem]{Lemma}

\newtheorem{question}{Question}

{
  \theoremstyle{definition}

}
\theoremstyle{remark}
\newtheorem{remark}{Remark}[section]

\numberwithin{equation}{section}

\begin{document}

\title{Fractal Sumset Properties}

\author{Derong Kong}
\address[D. Kong]{College of Mathematics and Statistics, Center of Mathematics, Chongqing University, Chongqing 401331, People's Republic of China}
\email{derongkong@126.com}

\author{Zhiqiang Wang*}
\address[Z. Wang]{School of Mathematical Sciences, Key Laboratory of MEA (Ministry of Education) \& Shanghai Key Laboratory of PMMP, East China Normal University, Shanghai 200241, People's Republic of China}
\email{zhiqiangwzy@163.com}
\date{\today}

\begin{abstract}
In this paper we introduce two notions of fractal sumset properties.
A compact set $K\subset\R^d$ is said to have the \emph{Hausdorff sumset property} (HSP) if for any $\ell\in\N_{\ge 2}$ there exist compact sets $K_1, K_2,\ldots, K_\ell$ such that $K_1+K_2+\cdots+K_\ell\subset K$ and $\dim_H K_i=\dim_H K$ for all $1\le i\le \ell$.
Analogously, if we replace the Hausdorff dimension by the packing dimension in the definition of HSP, then the compact set $K\subset\R^d$ is said to have the \emph{packing sumset property} (PSP).
We show that the HSP fails for certain homogeneous self-similar sets satisfying the strong separation condition, while the PSP holds for all homogeneous self-similar sets in $\R^d$.
\end{abstract}

\thanks{* Corresponding author}

\maketitle

\section{Introduction}
Let $\N$ be the set of natural numbers. The famous Erd{\H o}s sumset conjecture  states that if a set $A\subset\N$ has positive upper Banach density, then there exist two infinite sets $B, C\subset\N$ such that $A$ contains the \emph{sumset} $B+C:=\set{b+c: b\in B, c\in C}$ (see \cite{Erdos-Graham-1980}). This conjecture was fully proven by Moreira et al.~\cite{Moreira-Richter-Robertson-2019} in a more general setting including countable amenable groups. A short proof of this conjecture was later given by Host \cite{Host-2019}. Recently, Kra et al.~\cite{Kra-Moreira-Richter-Robertson-2022} extended this sumset result and showed that if $A\subset\N$ has positive upper Banach density, then for any $\ell\in\N_{\ge 2}$ there exist infinite sets $B_1, B_2,\ldots, B_\ell\subset \N$ such that $A$ contains  the sumset $B_1+B_2+\cdots+B_\ell:=\set{b_1+b_2+\cdots+b_\ell: b_i\in B_i,~1\le i\le \ell}$, where $\N_{\ge n}:=\set{\ell\in\N: \ell\ge n}$ for any $n\in\N$.

In the literature there is a great interest in the study of sumsets such as iterated sumsets $nE:=\underbrace{E+\cdots +E}_n$ and inhomogeneous sumsets $E_1+E_2+\cdots +E_n$. The study of  sumsets is closely related to the famous Marstrand's Projection Theorem (cf.~\cite{Marstrand-1954}). Fraser et al.~\cite{Fraser-Howroyd-Yu-2019} showed that if $E\subset\R$ is a closed set with positive lower dimension (see \cite{Fraser-2021} for its definition) then $\dim_H nE$ tends to $1$ as $n\to\f$.
Lindenstrauss et al.~\cite{Lindenstrauss-Miri-Peres-1999} considered the dimension growth for the inhomogeneous sumsets $E_1+E_2+\cdots+E_n$, where each $E_i$ is a compact $\times p$ invariant subset of the unit circle with positive Hausdorff dimension.
Recently, Feng and Wu \cite{Feng-Wu-2021} studied for which set $F\subset \R^d$ its iterated sumset $nF$ has non-empty interior for sufficiently large $n\in\N$. For more information on sumsets we refer to the papers \cite{Orponen-2022, Rossi-Shmerkin-2020} and the references therein.

Inspired by the Erd{\H o}s sumset conjecture and recent progress, we consider a fractal analogues.
Given a compact set $K \subset \R^d$ with positive Hausdorff dimension, it is interesting to ask whether $K$ contains a sumset $B+C$, where $B$ and $C$ are large in the sense of cardinality or dimension?

There do exist full dimension sets in $\R$ that have no sumset properties.
Keleti constructed in \cite{Keleti-1998} a compact set $A\subset\R$ such that $\dim_H A = 1$ and the intersection $A\cap (A+t)$ contains at most one point for all $t \in \R \setminus \{0\}$, where $A+t:=A+\set{t}$. We claim that the set $A$ has no sumset properties. To see this, suppose on the contrary that $A\supset A_1+A_2$ for some nonempety sets $A_1, A_2\subset\R$ with $\# A_1,\# A_2 \ge 2$, where $\#$ denotes the cardinality of a set.
Then taking $t_1 \ne t_2 \in A_1$ we have $A_2 + t_1 \subset A \cap (A+t_1 - t_2)$, a contradiction.
Thus, the set $A$ has no sumset properties.

Given $\ell\in\N_{\ge 2}$, a compact set $K\subset\R^d$  is said to have the \emph{Hausdorff $\ell$-sumset property} (simply called, HSP-$\ell$), if there exist compact sets $K_1, K_2,\ldots, K_\ell\subset\R^d$ such that
\begin{equation}\label{eq:HSP}
  K_1 + K_2 + \cdots + K_\ell \subset K, \quad\text{and}\quad \dim_H K_i =\dim_H K\quad\forall 1\le i\le \ell.
\end{equation}
If the set $K$ has HSP-$\ell$ for all $\ell\in\N_{\ge 2}$, then we say that $K$ has the \emph{Hausdorff sumset property} (shortly called, HSP).
Analogously, a compact set $K\subset\R^d$ is said to have the \emph{packing $\ell$-sumset property} (simply called, PSP-$\ell$), if there exist compact sets $K_1, K_2, \ldots, K_\ell\subset\R^d$ satisfying (\ref{eq:HSP}) with the Hausdorff dimension replaced by the packing dimension. Similarly, a set $K\subset\R^d$ is said to have the \emph{packing sumset property} (shortly called, PSP) if $K$ has the PSP-$\ell$ for all $\ell\in\N_{\ge 2}$.

A classical result by Erd{\H o}s and Volkmann \cite{Erdos-Volkmann-1966} showed that for any $\alpha\in(0,1)$ there exists an additive subgroup $G$ of $\R$ with $\dim_H G=\alpha$.
By modifying Erd{\H o}s--Volkmann's example, we can construct compact sets with any given Hausdorff dimension having both the HSP and the PSP.
Fix $\alpha \in(0, 1)$.
For $k \ge 1$, let $m_k = (k+1)!$  and   $N_k = \lfloor m_k^{1/\alpha} \rfloor$, where $\lfloor x \rfloor$ denotes the integer part of $x$.
Define
\begin{equation}\label{set-K}
  K = \bigg\{ \sum_{k=1}^{\f} \frac{d_k}{N_1 N_2 \cdots N_k}:\; d_k \in \{0,1,\ldots, m_k-1\} \; \forall k \ge 1 \bigg\}.
\end{equation}
For $\ell \in \N_{\ge 2}$, let $m_k' := m_k / \ell=(k+1)!/\ell$ for $k > \ell$, and set
\begin{equation*}
  B_\ell = \bigg\{ \sum_{k=\ell+1}^{\f} \frac{d_k}{N_1 N_2 \cdots N_k}:\; d_k \in \{0,1,\ldots, m'_k-1\} \; \forall k > \ell \bigg\}.
\end{equation*}
It is straightforward to check the iterated sumset
$$ \ell B_\ell = \bigg\{ \sum_{k={\ell+1}}^{\f} \frac{d_k}{N_1 N_2 \cdots N_k}:\; d_k \in \{0,1,\ldots, \ell({m'_k}-1) \} \; \forall k{>} \ell \bigg\} \subset K. $$
Furthermore, we have $\dim_H K = \dim_P K = \dim_H B_\ell = \dim_P B_\ell =\alpha$.
The dimensions of   $K$ and $B_\ell$ can be explicitly calculated by the dimension formulae of partial homogeneous Cantor sets due to Feng et al.~\cite{Feng-Wen-Wu-1997}.
Thus, the set $K$ defined in (\ref{set-K}) has both the HSP and the PSP.
This example can be easily extended to higher dimensions   by taking the Cartesian products.

In this paper, we focus on the sumset properties in self-similar sets.
A non-empty compact set $K \subset \R^d$ is called a \emph{self-similar set} if there exists an \emph{iterated function system} (IFS) $\Phi=\{ f_j(x)=\rho_j O_j x+b_j \}_{j=1}^m$ with each $\rho_j\in(0,1)$, $b_j \in \R^d$, and $O_j$ a $d\times d$ orthogonal real matrix, such that $K=\bigcup_{j=1}^m f_j(K)$ (cf.~\cite{Falconer-2014,Hutchinson_1981}).
If all $\rho_j O_j, 1\le j\le m,$ are identical, then we say $K$ is a \emph{homogeneous} self-similar set.
We say that the IFS $\Phi$ satisfies the \emph{strong separation condition} (SSC) if the union $\bigcup_{j=1}^m f_j(K)$ is pairwise disjoint.
To avoid the trivial case, we always assume the self-similar set not to be a singleton.

Our first result states that a class of homogeneous self-similar sets in $\R^d$ does not have the HSP-$2$, and hence fails the HSP.
Let $K \subset \R^d$ be a homogeneous self-similar set generated by the IFS $\Phi = \{f_b(x) = \rho O x + b: b \in D\}$, where $D \subset \R^d$ is a finite set with $\# D \ge 2$.
Then the set $K$ can be written as $$K = \bigg\{ \sum_{k=1}^{\f} (\rho O)^{k-1} b_k: b_k \in D \;\forall k \ge 1 \bigg\}.$$
Whence,$$K-K = \bigg\{ \sum_{k=1}^{\f} (\rho O)^{k-1} t_k: t_k \in D-D \;\forall k \ge 1 \bigg\}.$$
So, $K-K$ is a self-similar set generated by the IFS $\Psi = \{ g_b(x) = \rho O x + b: b \in D-D \}$.

\begin{theorem}\label{th:HSP-fail-self-similar}
  If the IFS $\Psi$ satisfies the SSC, then there exists a constant $\beta < \dim_H K$ such that for any two non-empty sets $K_1,K_2 \subset \R^d$ satisfying $K_1 + K_2 \subset K$ we have
  \[ \dim_H K_1 + \dim_H K_2 \le 2 \beta. \]
  Consequently, the set $K$ does not have the HSP-$2$, and hence fails the HSP.
\end{theorem}

\begin{remark}\label{rem:2}
  Note that Theorem \ref{th:HSP-fail-self-similar} can be strengthened in $\R$.
  By \cite[Theorem 2.1]{Fraser-Howroyd-Yu-2019}, for any set $K \subset \R$ with $0 < \dim_H K = \ol{\dim}_B K \le \dim_A K < 1$, if $K_1 + K_2 \subset K$ and $\dim_H K_1 = \dim_H K$ then we have $\dim_H K_2 = 0$.
  For the definition of upper Box dimension $\overline{\dim}_B$ and Assouad dimension $\dim_A$, we refer to the book of Fraser \cite{Fraser-2021}.
  So, \cite[Theorem 2.1]{Fraser-Howroyd-Yu-2019} implies that any self-similar set in $\R$ with non-integer dimension does not have the HSP-$2$, and hence fails the HSP.
\end{remark}

In contrast with Theorem \ref{th:HSP-fail-self-similar}, the PSP holds for all homogeneous self-similar sets, even without the SSC.
\begin{theorem}\label{th:PSP-self-similar}
  If $K\subset\R^d$ is a homogeneous self-similar set, then $K$ has the PSP, i.e., for any $\ell\in\N_{\ge 2}$ there exist compact subsets $K_1, K_2,\ldots, K_\ell\subset\R^d$ such that
  \begin{equation*}\label{eq:PSP}
  K_1 + K_2 + \cdots + K_\ell \subset K, \quad\text{and}\quad\dim_P K_i = \dim_P K \quad \forall 1\le i\le \ell.
  \end{equation*}
\end{theorem}

Another analogue is to ask whether a set $K\subset\R^d$ with positive dimension contains a sumset $B+C$ with both sets $B$ and $C$ having positive dimension.
A compact set $K\subset\R^d$ is said to have the \emph{positive dimension sumset property} (simply called, PDSP) if $\dim_H K>0$, and for any $\ell\in\N_{\ge 2}$ there exist compact sets $K_1, K_2,\ldots, K_\ell\subset\R^d$ such that
$$ K_1+K_2+\cdots+K_\ell\subset K, \quad\textrm{and}\quad \dim_H K_i>0\quad\forall 1\le i\le \ell. $$

Our final result shows that the PDSP holds for homogeneous self-similar sets.

\begin{theorem}\label{th:PDSP}
  If $K$ is a homogeneous self-similar set in $\R^d$, or a self-similar set in $\R$ or $\R^2$,
 then $K$ has the PDSP.
\end{theorem}

The rest of the paper is organized as follows.
In Section \ref{sec:HSP-fails} we prove that the HSP fails for a class of self-similar sets (Theorem \ref{th:HSP-fail-self-similar}). On the other hand, the PSP always holds for homogeneous self-similar sets (Theorem \ref{th:PSP-self-similar}). This will be proven in  Section \ref{sec:PSP-holds}.
Finally, in Section \ref{sec:PDSP} we consider the PDSP and prove Theorem \ref{th:PDSP}. We also propose some questions.

\section{The HSP is not universal}\label{sec:HSP-fails}
In this section we will show that the HSP is not universal, and prove Theorem \ref{th:HSP-fail-self-similar}.
First, we need the following combinatorial lemma.

\begin{lemma}\label{lem:digit-translation}
  Let $A \subset \R^d$ be a finite subset with $\# A \ge 2$.
  Then for any distinct points $t_1,t_2,\dots,t_\ell \in \R^d$ with $1 \le \ell \le \#A +1$, we have
  \[ \# \bigg( \bigcap_{j=1}^\ell (A + t_j) \bigg) \le \#A +1 - \ell. \]
\end{lemma}
\begin{proof}
  The inequality is clear for $\ell=1$. In the following we assume $\ell \ge 2$. We first define a lexicographical ordering for points in $\R^d$. For $x=(x_1,\dots,x_d),y=(y_1,\dots,y_d) \in \R^d$, we define $x \prec y$ if $x_1 < y_1$, or there exists $1\le j \le d-1$ such that $x_1= y_1$, $\ldots$, $x_j = y_j$, and $x_{j+1} < y_{j+1}$.
  It is easy to verify that if $x\prec y$ and $x' \prec y'$ then we have $x+x' \prec y+y'$.

  We write $m:=\# A$ and $A=\{a_1, a_2, \dots, a_m\}\subset\R^d$ with $a_1 \prec a_2 \prec \cdots \prec a_m$.
  Without loss of generality, we can assume that $t_1 \prec t_2 \prec \cdots \prec t_\ell$.
  It suffices to prove that
  \[\set{t_1 + a_1, \dots, t_1 + a_{\ell-1}}\cap\bigcap_{j=1}^\ell(A+t_j)=\emptyset.\]
  Suppose on the contrary that $t_1 + a_k\in\bigcap_{j=1}^\ell(A+t_j)$ for some $1\le k \le \ell-1$.
  Then there exist $k_2, \dots, k_\ell \in \{1,2,\dots,m\}$ such that
  \[ t_1 + a_k = t_2 + a_{k_2} = \cdots = t_\ell + a_{k_\ell}. \]
  Since $t_1 \prec t_2 \prec \cdots \prec t_\ell$, we have $a_{k_\ell} \prec \cdots\prec a_{k_2} \prec a_k$.
  It follows that $k \ge \ell$, a contradiction.
\end{proof}

Recall that $K \subset \R^d$ is a homogeneous self-similar set generated by the IFS $\Phi = \{f_b(x) = \rho O x + b: b \in D\},$ where $D \subset \R^d$ is a finite set with $\# D \ge 2$.
Then the set $K$ can be written as \[ K = \bigg\{ \sum_{k=1}^{\f} (\rho O)^{k-1} b_k: b_k \in D \; \forall k \ge 1 \bigg\}. \]
The difference set
\[ K-K = \bigg\{ \sum_{k=1}^{\f} (\rho O)^{k-1} t_k: t_k \in D-D \; \forall k \ge 1 \bigg\} \]
is also a self-similar set generated by the IFS $\Psi = \{ g_b(x) = \rho O x + b: b \in D-D \}$.
Suppose the IFS $\Psi$ satisfies the SSC.
Then for each $t\in K-K$ there exists a unique sequence $(t_k) \in (D-D)^\N$ such that \[ t = \sum_{k=1}^{\f} (\rho O)^{k-1} t_k, \]
and the unique sequence $(t_k)$ is called a coding of $t$ with respect to the digit set $D-D$. In this case, we have
\begin{equation}\label{eq:translation-intersection}
  K \cap (K+t) = \bigg\{ \sum_{k=1}^{\f} (\rho O)^{k-1} b_k: b_k \in D\cap(D+t_k) \; \forall k \ge 1 \bigg\}.
\end{equation}

\begin{proof}[Proof of Theorem \ref{th:HSP-fail-self-similar}]
  Without loss of generality, we assume that $\mathbf{0}=(0,0,\ldots, 0) \in D$. {Let $K_1$ and $K_2$ be two non-empty subsets of $\R^d$ satisfying $K_1+K_2\subset K$.}

  Take $x_0 \in K_1$, and then we have \[ (K_1 - x_0) + (K_2 + x_0) \subset K. \]
  Note that $\mathbf{0} \in K_1 - x_0$, and Hausdorff dimension is stable under translations.
  Then we can assume without loss of generality that $\mathbf{0} \in K_1$.
  So, by using $K_1+K_2\subset K$ it follows that $K_2 \subset K$, and hence,
  \begin{equation}\label{eq:K-1-K-2}
     K_1\subset K-K,\qquad
     K_2 \subset  \bigcap_{t \in -K_1} \big( K \cap (K+t) \big).
  \end{equation}

  Note that $t \in -K_1 \subset K - K$. Since the IFS $\Psi$ satisfies the SSC, each $t\in -K_1$ has a unique $D-D$ coding.
  Let $\Lambda$ denote the set of all unique codings $(t_k) \in (D-D)^\N$ of points in $-K_1$.
  For $k \ge 1$, let $\Lambda_k \subset D-D$ be the set of all possible digits occurring in the $k$-th position of sequences in $\Lambda$. Then $\Lambda\subset\prod_{k=1}^\f\Lambda_k$.
  Since $\mathbf{0} \in K_1$ and the IFS $\Psi$ satisfies the SSC, we have $0 \in \Lambda_k$ for all $k \ge 1$.
  By (\ref{eq:translation-intersection}) and (\ref{eq:K-1-K-2}), we have
  \[ K_2 \subset \wt{K_2}: = \bigg\{ \sum_{k=1}^{\f} (\rho O)^{k-1} b_k: b_k \in \bigcap_{b \in \Lambda_k}(D+b) \;\; \forall k \ge 1 \bigg\}. \]
  For $k \ge 1$, we write \[ m_k = \# \bigg( \bigcap_{b \in \Lambda_k}(D+b) \bigg). \]
  Then we obtain \begin{equation}\label{eq:upper-K2} \dim_H K_2 \le \dim_H \wt{K_2} \le \liminf_{k \to \f} \frac{\log (m_1 m_2 \cdots m_k)}{-k\log \rho}=\liminf_{k\to\f}\frac{\sum_{j=1}^k\log m_j}{-k\log\rho}.  \end{equation}

  On the other hand, since each point in $-K_1$ has a unique coding in  $\Lambda \subset \prod_{k=1}^\f\Lambda_k$, we have
  \begin{equation}\label{eq:upper-K1} \dim_H K_1 = \dim_H (-K_1) \le \liminf_{k \to \f} \frac{\sum_{j=1}^k\log (\#\Lambda_j)}{-k\log\rho}. \end{equation}
  Note by Lemma \ref{lem:digit-translation} that $m_j+ \# \Lambda_j\le \# D + 1$ for all $j \ge 1$. Therefore, by (\ref{eq:upper-K2}) and (\ref{eq:upper-K1}) it follows that
  \begin{align*}
    \dim_H K_1 + \dim_H  K_2
    & \le \liminf_{k\to\f}\frac{\sum_{j=1}^k(\log m_j+\log(\#\Lambda_j))}{-k\log\rho} \\
    & \le \liminf_{k\to\f}\frac{\sum_{j=1}^k(\log m_j+\log(\# D+1-m_j))}{-k\log\rho} \\
    & \le \frac{\gamma}{-\log \rho},
  \end{align*}
  where
  \begin{equation}\label{eq:def-gamma}
    \gamma:= \max\big\{ \log m + \log(\#D +1 - m) : m \in \{1,2,\dots,\#D\} \big\}.
  \end{equation}
  Observe by the concavity of the function $\log x$ that
  \[ \gamma \le 2 \log \left(\frac{\# D + 1}{2}\right) < 2 \log (\# D), \]
  where the second inequality follows by $\#D\ge 2$.
  Hence, \[ \dim_H K_1 + \dim_H  K_2 \le \frac{\gamma}{-\log\rho}<2\frac{\log\#D}{-\log\rho}=2\dim_H K.\]
  This completes the proof by setting $\beta = \frac{ \gamma}{-2 \log \rho}$.
\end{proof}

\begin{remark}
  When $\# D = 2$, the number $\gamma$ defined in (\ref{eq:def-gamma}) is indeed $\log \#D$, and then we can conclude that \[ \dim_H K_1 + \dim_H K_2 \le \dim_H K \] for any two non-empty subsets $K_1,K_2 \subset \R^d$ satisfying $K_1 + K_2 \subset K$.
\end{remark}

At the end of this section we point out that Theorem \ref{th:HSP-fail-self-similar} can be applied to  homogeneous self-similar sets in $\R$.
For a positive integer $N$ and a real number $0< \rho < 1/(N+1)$, we write
\[ E_{\rho,N} = \bigg\{ \frac{1-\rho}{N} \sum_{k=1}^{\f} x_k \rho^{k-1}: x_k \in \{0,1,\ldots, N\}\; \forall k \ge 1 \bigg\}. \]
By Theorem \ref{th:HSP-fail-self-similar},  for $0 < \rho < 1/(2N+1)$ the set $E_{\rho, N}$ does not have the HSP-$2$, and hence fails the HSP.

\section{Homogeneous self-similar set has the PSP}\label{sec:PSP-holds}
In contrast with Theorem \ref{th:HSP-fail-self-similar} we show in this section that the PSP always holds for homogeneous self-similar sets, and prove Theorem \ref{th:PSP-self-similar}.
Recall that $K$ is a homogeneous self-similar set in $\R^d$ if it can be generated by an IFS $\{ f_j(x) = \rho O x + b_j \}_{j=1}^m$, where $0< \rho <1$, $O$ is a $d \times d$ orthogonal real matrix, and each $b_j\in\R^d$. Without loss of generality  we may assume that $b_1 = {\bf 0}$.

Let $\Omega = \{ 1,2,\ldots, m\}^\N$ be the set of all infinite sequences $(i_k)$ with each digit $i_k\in\set{1,2,\ldots, m}$. Equipped with the product topology of the discrete topology on $\set{1,2,\ldots, m}$, $\Omega$ becomes a compact metric space.
We define the coding map $\pi: \Omega \to K$  by
\begin{align*}
  \pi((i_k)) & = \lim_{k \to \f} f_{i_1} \circ f_{i_2} \circ \cdots \circ f_{i_k} ({\bf 0})  = \sum_{k=1}^{\f} (\rho O)^{k-1} b_{i_k}.
\end{align*}
Then $\pi$ is continuous and surjective.
For $S \subset \N$, we define \[ \Omega_S := \big\{ (i_k) \in \Omega: \;  i_k =1 \text{ for } k \notin S \big\}\quad \text{ and }\quad K_S := \pi\big( \Omega_S \big). \]
Then $K_S$ is a subset of $K=\pi(\Omega)$.

\begin{proposition}\label{prop:packing-dimension}
  If $S \subset \N$ satisfies
  \[ \limsup_{n \to \f} \frac{\#(S \cap [1,n])}{n} =1, \]
  then we have \[ \dim_P K_S=\overline{\dim}_B K_S = \dim_P K. \]
\end{proposition}
\begin{proof}
  For the first equality, let $V$ be an open subset that intersects $K_S$.
  Then we can find $(j_n) \in \Omega_S$ such that $ \pi( (j_n) ) \in V. $
  Since $V$ is open and $\pi$ is continuous, there exists $n_0\in\N$ such that
  $\pi( [j_1 j_2 \cdots j_{n_0}] ) \subset V, $
  where $[j_1 j_2 \cdots j_{n_0}] := \{ (i_n) \in \Omega: i_k = j_k \text{ for } 1 \le k \le n_0 \}$ is a cylinder set.
  It follows that \[ \pi\big( [j_1 j_2 \cdots j_{n_0}] \cap \Omega_S \big) \subset  V \cap K_S. \]
  Note that the set $K_S$ is the union of finitely many translations of   $\pi\big( [j_1 j_2 \cdots j_{n_0}] \cap \Omega_S \big)$, and the upper box-counting dimension is finitely stable.
  Therefore, we conclude that \[ \ol{\dim}_B (V \cap K_S) = \ol{\dim}_B K_S. \]
  Note that $K_S$ is compact.
  Thus, by \cite[Corollary 3.10]{Falconer-2014} (see also \cite[Corollary 2.8.2]{Bishop-Peres-2017}) it follows that \[ \dim_P K_S = \ol{\dim}_B K_S \]as desired.

  For the second equality, we first observe that $\ol{\dim}_B K_S \le \ol{\dim}_B K = \dim_P K$. So it suffices to prove the inverse inequality. Without loss of generality we assume that $\diam K=1$.
  Write $\alpha := \dim_P K$.
  Suppose {on the contrary} that $ \overline{\dim}_B K_S  < \alpha. $
  Then there exists $\ep_0\in(0, \alpha)$ such that for all sufficiently large $n$, we have \begin{equation}\label{eq:august-9-1} \frac{\log N_{\rho^n}(K_S)}{-n\log \rho} \le \alpha -\ep_0, \end{equation}
  where $N_\delta(F)$ denotes the smallest number of closed balls of radius $\delta$ that cover a set $F\subset \R^d$.

  For $n \ge 1$, define \[ \Omega_{S,n} := \big\{ (i_k) \in \Omega: \;  i_k =1 \text{ for } k \notin S\cup \N_{\ge n+1} \big\}\quad \text{ and }\quad K_{S,n} := \pi\big( \Omega_{S,n} \big). \]
  That is, $\Omega_{S,n}$ is the union of all $n$-level cylinder sets that intersects $\Omega_S$.
  Note by $\diam K=1$ that the diameter of all $n$-level sets of $K$ is $\rho^n$.
  Thus we have \[ N_{3\rho^n}(K_{S,n}) \le N_{\rho^n}(K_S). \]
  Note that the set $K$ is covered by the union of at most $ m^{n- \#(S \cap [1,n])} $
  many translations of $K_{S,n}$.
  It follows that \begin{equation}\label{eq:august-9-2} N_{3\rho^n}(K)\le m^{n- \#(S \cap [1,n])} \cdot N_{3\rho^n}(K_{S,n}) \le m^{n- \#(S \cap [1,n])} \cdot N_{\rho^n}(K_S). \end{equation}
  Thus, {by (\ref{eq:august-9-1}) and (\ref{eq:august-9-2})} we have
  \begin{align*}
    \ul{\dim}_B K & \le \liminf_{n \to \f} \frac{\log N_{3\rho^n}(K)}{-\log(3\rho^n)} \\
    & \le \liminf_{n \to \f} \frac{ \log \Big( m^{n- \#(S \cap [1,n])} \cdot N_{\rho^n}(K_S) \Big) }{- n \log \rho} \\
    & \le \alpha - \ep_0 + \frac{\log m}{-\log \rho} \cdot \liminf_{n \to \f} \Big( 1 - \frac{\#(S \cap [1,n])}{n} \Big) \\
    & = \alpha-\ep_0.
  \end{align*}
  This contradicts the fact that $\ul{\dim}_B K = \ol{\dim}_B K =\dim_P K = \alpha$ by the self-similarity of $K$ (cf. \cite[Corollary 3.3]{Falconer-1997}).
  Thus, we obtain $\ol{\dim}_B K_S \ge \alpha$ as desired.
\end{proof}

\begin{lemma}\label{lem:divide-limsup}
  Let $S \subset \N$ satisfy \begin{equation}\label{eq:sup-assum} \beta = \limsup_{n \to \f} \frac{\#(S \cap [1,n])}{n} >0. \end{equation}
  Then for any $\ell \in \N_{\ge 2}$ we can divide $S$ into {pairwise} disjoint subsets $S_1, S_2, \ldots, S_\ell$ such that
  \[ \limsup_{n \to \f} \frac{\#(S_j \cap [1,n])}{n} = \beta \quad \forall 1 \le j \le \ell. \]
\end{lemma}
\begin{proof}
  Note that for any given $m \ge 0$, we   have \[ \limsup_{n \to \f} \frac{\#(S \cap [m + 1,n])}{n} = \beta. \]
  Thus, we can define an increasing sequence $\{n_k\}$ of integers such that $n_0 = 0$ and
  \begin{equation}\label{eq:august-9-3}
    \frac{\#(S \cap [n_{k-1}+1, n_k])}{n_k} > \beta - \frac{1}{k} \quad \forall k \ge 1.
  \end{equation}

  Fix $\ell \in \N_{\ge 2}$.
  For $1 \le j \le \ell$, we define
  \[ S_j = \bigcup_{k=0}^\f \big( S \cap [n_{k\ell+j-1} +1,  n_{k\ell+j}] \big). \]
  Clearly, $S_1, S_2, \ldots, S_\ell$ are pairwise disjoint, and $ S = \bigcup_{j=1}^\ell S_j. $
  Note {by (\ref{eq:august-9-3})} that for each $j\in\set{1,2,\ldots, \ell}$, \[ \frac{\#(S_j \cap [1,n_{k\ell+j}])}{n_{k\ell+j}} > \beta - \frac{1}{k\ell +j}\quad\forall k\in\N. \]
  From this and (\ref{eq:sup-assum}) we conclude that \[ \limsup_{n \to \f} \frac{\#(S_j \cap [1,n])}{n} = \beta \quad \forall 1 \le j \le \ell. \]
\end{proof}

\begin{proof}[Proof of Theorem \ref{th:PSP-self-similar}]
  Fix $\ell \in \N_{\ge 2}$.
  By Lemma \ref{lem:divide-limsup}, we can divide $\N$ into pairwise disjoint subsets $S_1, S_2, \ldots, S_\ell$ such that
  \[ \limsup_{n \to \f} \frac{\#(S_j \cap [1,n])}{n} = 1 \quad \forall 1 \le j \le \ell. \] By Proposition \ref{prop:packing-dimension} this implies that \[ \dim_P K_{S_j} =  \dim_P K \quad \forall 1 \le j \le \ell, \]
  Note that $b_1={\bf 0}$. By using the pairwise disjointness of $S_1, S_2, \ldots, S_\ell$ it is easy to check that \[ K_{S_1} + K_{S_2} + \cdots + K_{S_\ell} =K. \]
  Thus, the set $K$ has the PSP-$\ell$.
  Since $\ell$ was arbitrary, we conclude that the homogeneous self-similar $K$ has the PSP.
\end{proof}

\section{The PDSP and questions}\label{sec:PDSP}

In this section we will consider the PDSP and prove Theorem \ref{th:PDSP}.
The following lemma   can be deduced from  \cite[Corollary 1.1]{Xie-Yin-Sun-2003}.
\begin{lemma}
  \label{lem:positive-dim}
  If $K\subset\R^d$ is a self-similar set containing at least two points, then $\dim_H K>0$.
\end{lemma}

First, we show that any homogeneous self-similar set has the PDSP.
As a consequence, any compact set containing a homogeneous self-similar set has the PDSP.

\begin{lemma}\label{lemma:homogeneous-PDSP}
If $K$ is a homogeneous self-similar set in $\R^d$, then  $K$ has the PDSP.
\end{lemma}
\begin{proof}
 Let $K$ be a homogeneous self-similar set generated by $\{ f_j(x) = \rho O x + b_j \}_{j=1}^m$. Note by our assumption that $K$ contains at least two points.
  Then by Lemma \ref{lem:positive-dim} we have $\dim_H K >0$.
  Take $\ell\in\N_{\ge 2}$. By choosing $S_j=\ell\N+j-1$ with $j\in\{1,2,\ldots,\ell\}$ in the proof of Theorem \ref{th:PSP-self-similar} one can verify that $K_{S_1}+K_{S_2}+\cdots + K_{S_\ell}\subset K$.
  Furthermore, since $S_1= \ell \N$, the set $K_{S_1}$ is the self-similar set generated by the IFS $$\big\{ g_j(x) = (\rho O)^\ell x + b_j \big\}_{j=1}^m, $$
  which implies $\dim_H K_{S_1}>0$ by Lemma \ref{lem:positive-dim}.
  Note that $K_{S_j} = (\rho O)^{j - 1} K_{S_1}$ for all $1 < j \le \ell$.
  Thus, $$\dim_H K_{S_1} = \dim_H K_{S_2} = \cdots = \dim_H K_{S_\ell} > 0.$$
  Therefore, $K$ has the PDSP.
\end{proof}

\begin{proof}[Proof of Theorem \ref{th:PDSP}]
  By Lemma \ref{lemma:homogeneous-PDSP}, it remains to show that any self-similar set in $\R$ or $\R^2$ contains a homogeneous self-similar set.
  Let $K$ be a self-similar set in $\R$ or $\R^2$ generated by the IFS $\{ f_j(x) = \rho_j O_j x + b_j\}_{j=1}^m$.
  Without loss of generality, we assume that $f_1$ and $f_2$ have distinct fixed points.

  In $\R$, the homogeneous self-similar set generated by $\{ f_1 \circ f_2, f_2 \circ f_1 \}$ is contained in $K$.
  In $\R^2$, all the rotation matrices are commutative.
  Note that   $f_1\circ f_1$ and $f_2\circ f_2$ contain no reflections, and their   fixed points are distinct.
  Thus the homogeneous self-similar set generated by $$\{ f_1\circ f_1 \circ f_2\circ f_2, \; f_2\circ f_2\circ f_1\circ f_1 \}$$ is contained in $K$.
  This completes the proof.
\end{proof}

It is worth mentioning that the PDSP may hold for a even larger class of sets, for example the set of all graph-directed sets in $\R$ or $\R^2$ (cf.~\cite{Farkas-2019}).

\medskip

At the end of this section we propose some questions. Note by  Theorem \ref{th:HSP-fail-self-similar} that the HSP fails for a class of homogeneous self-similar set satisfying the SSC.
\begin{question}
  Do all self-similar sets in $\R^d$ with non-integer dimension fail the HSP?
\end{question}

In view of Theorem \ref{th:PDSP}, the PDSP holds for all self-similar sets in $\R$ or $\R^2$. Furthermore, it holds for homogeneous self-similar sets in $\R^d$ with $d\ge 3$.
\begin{question}
  Do  all inhomogeneous self-similar sets in $\R^d, d \ge 3,$ have the PDSP?
\end{question}

Note by Lebesgue density theorem that a Lebesgue measurable set of positive measure contains a similarity copy of any finite set. It might be interesting to ask the following question.

\begin{question}
  Given a compact set $K \subset \R^d$ with positive Lebesgue measure, do there exist two infinite sets $B, C \subset \R^d$ such that $B+C \subset K?$
\end{question}

\section*{Acknowledgements}
The authors thank the anonymous referees for many useful suggestions which improve the presentation of the paper.
The first author wants to thank Tuomas Orpenon for many useful discussions during the Fractal Geometry conference in Edinburgh 2023, especially for Remark \ref{rem:2} and many useful references.
The first author was supported by NSFC No.~11971079. The second author was supported by NSFC No.~12071148, and Science and Technology Commission of Shanghai Municipality (STCSM) No.~22DZ2229014, and Fundamental Research Funds for the Central Universities No.~YBNLTS2023-016.

\end{document}